\documentclass[]{interact}

\usepackage{epstopdf}
\usepackage[caption=false]{subfig}
\usepackage{multirow}

\usepackage[numbers,sort&compress]{natbib}
\bibpunct[, ]{[}{]}{,}{n}{,}{,}
\renewcommand\bibfont{\fontsize{10}{12}\selectfont}

\theoremstyle{plain}
\newtheorem{theorem}{Theorem}[section]
\newtheorem{lemma}[theorem]{Lemma}
\newtheorem{corollary}[theorem]{Corollary}

\theoremstyle{definition}
\newtheorem{definition}[theorem]{Definition}

\theoremstyle{remark}

\begin{document}

\articletype{Manuscript}

\title{A feasible adaptive refinement algorithm for linear semi-infinite optimization}

\author{
\name{Shuxiong Wang\textsuperscript{1}\thanks{Email: shuxionw@uci.edu}}
\affil{\textsuperscript{1} Department of Mathematics, University of California, Irvine, CA,US}
}

\maketitle

\begin{abstract}
A numerical method is developed to solve linear semi-infinite programming problem (LSIP) in which
the iterates produced by the algorithm are feasible for the original problem.
This is achieved by constructing a sequence of standard linear programming problems 
with respect to the successive discretization of the index set such that the approximate
regions are included in the original feasible region. 
The convergence of the approximate solutions to the solution of the original problem is proved 
and the associated optimal objective function values of the approximate problems are 
monotonically decreasing and converge to the optimal value of LSIP.
An adaptive refinement procedure is designed to discretize the index set and update the constraints
for the approximate problem. 
Numerical experiments demonstrate the performance of the proposed algorithm.
\end{abstract}


\section{Introduction}
Linear semi-infinite programming problem (LSIP) refers to the optimization problem with finitely many
decision variables and infinitely many linear constraints associated with some parameters, 
which can be formulated as
\begin{gather}\label{LSIP}
 \begin{split}
   \min\limits_{x\in \mathbb{R}^n}\quad &  c^{\top}x                  \\
    \textrm{s.t.}\quad                  & a(y)^{\top}x + a_0(y)\geq 0\ \forall y\in Y, \\
                                        & x_i \geq 0, i = 1,2,...,n,
 \end{split}\tag*{(LSIP)}
\end{gather}
where $c\in \mathbb{R}^n$, $a(y) = [a_1(y),...,a_n(y)]^{\top}$ and 
$a_i: \mathbb{R}^m \mapsto \mathbb{R}$, for $i = 0,1,...,n$,
are real-valued coefficient functions, $Y \subseteq \mathbb{R}^m$ is the index set.
In this paper, we assume that $Y = [a,b]$ is an interval with $a<b$.
Denote by $F$ the feasible set of (\ref{LSIP}):
\begin{equation*}
  F = \{x\in \mathbb{R}^n_{+}\ |\ a(y)^{\top}x + a_0(y)\geq 0, \forall y\in Y\},
\end{equation*}
where $\mathbb{R}^n_{+} = \{x\in \mathbb{R}^n\ |\ x_i \geq 0, i = 1,2,...,n\}$.

Linear semi-infinite programming has wide applications in economics, robust optimization and numerous
engineering problems, etc.
More details can be found in \cite{yy1,yy2,yy3} and references therein.

Numerical methods have been proposed for solving linear semi-infinite programming problems
such as discretization methods, local reduction methods and descent direction methods
 (See \cite{mma,ma,hkko,rg} for an overview of these methods).
The main idea of discretization methods is to solve the following
linear program 
\begin{align*}
  \min_{x\in \mathbb{R}^n_{+}}  &\quad f(x) \\
  \textrm{s.t.}             &\quad a(y)^{\top}x + a_0(y)\geq 0\ \forall y\in T,
\end{align*}
in which the original index set $Y$ in \ref{LSIP} is replaced by its finite subset $T$.
The iterates generated by the discretization methods converge to a solution of the original
problem as the distance between $T$ and $Y$ tends to zero (see \cite{BB,yy2,mma}).
The reduction methods solve nonlinear equations by qusi-Newton method,
which require the smoothing conditions on the functions defining the constraint \cite{sa}.
The feasible descent direction methods generate a feasible direction based on the
current iterate and achieve the next iterate by such a direction \cite{tse}.

The purification methods proposed in \cite{ea,em} generate a finite feasible sequence where the
objective function value of each iterate is reduced.
The method proposed in \cite{ea} requires that the feasible set of \ref{LSIP} is locally polyhedral,
and the method proposed in \cite{em}
requires that the coefficient functions $a_i, i=0,1,...,n,$ are analytic.

Feasible iterative methods for nonlinear semi-infinite optimization problems have been developed
via techniques of convexification or concavification etc \cite{floudas,wang,Amit}. 
These methods might be applicable to solve \ref{LSIP} directly. However,  they are not
developed specifically for \ref{LSIP}. Computational time will be reduced if the algorithm
can be adapted to linear case effectively. 

In this paper, we develop a feasible iterative algorithm to solve \ref{LSIP}. The basic idea is to 
construct a sequence of standard linear optimization problems with respect to
the discretized subsets of the index set such that the feasible region of each linear optimization
problem is included in the feasible region of \ref{LSIP}. The
proposed method consists of two stages. The first stage is based on the restriction of the
semi-infinite constraint.
The second  stage is base on estimating the lower bound of the coefficient functions using 
concavification or interval method.

The rest of the paper is organized as follows. In section 2, we propose the methods to construct
the inner approximate regions for the feasible region of \ref{LSIP}. Numerical method to solve the
original linear semi-infinite programming problem is proposed in section 3. In section 4, we
implement our algorithm to some numerical examples to show the performance of the method. At last,
we conclude our paper in section 5.

\section{Restriction of the lower level problem}
The restriction of the lower level problem leads to inner approximation of the feasible region of
\ref{LSIP}, and thus, to feasible iterates. Two-stage procedures are performed to achieve the restriction for \ref{LSIP}. 
At the first stage, we construct an uniform lower-bound function w.r.t decision variables 
for the function defining constraint in \ref{LSIP}. This step requires to solve a global optimization
associated with coefficient functions over the index set. 
The second stage is to estimate the lower bound of the coefficient functions over the index set rather
than solving the optimization problems globally which significantly reduce the computational cost. 

\subsection{Construction of the lower-bound function}
The semi-infinite constraint of \ref{LSIP} can be reformulated as
\begin{equation}\label{fs1}
 \min_{y\in Y}  \{a(y)^{\top} x + a_0(y)\} \geq 0.
\end{equation}
Since $a(y)^{\top} x = \sum_{i=1}^{n} a_i(y) x_i$, (\ref{fs1}) is equivalent to
\[\min_{y\in Y}\{ \sum_{i=1}^{n} a_i(y) x_i + a_0(y)\}\geq 0.\]
By exchanging the minimization and summation on the left side of the inequality, we 
obtain a new linear inequality
\begin{equation}
  \sum_{i=1}^{n} \{\min_{y\in Y} a_i(y)\} x_i + \min_{y\in Y} a_0(y)\geq 0.
  \label{fs2}
\end{equation}
Since the decision variables $x_i \geq 0$, $i = 1,2,...,n$, we have
\[\sum_{i=1}^{n} \{\min_{y\in Y} a_i(y)\} x_i + \min_{y\in Y} a_0(y) \leq \min_{y\in Y}\{ \sum_{i=1}^{n} a_i(y) x_i + a_0(y)\}.\]
Thus, we obtain an uniform lower-bound function for $\min_{y\in Y}  \{a(y)^{\top} x + a_0(y)\} $. 
And any point $x$ satisfying (\ref{fs2}) is a feasible point for LSIP.
Let $\bar{F}$ be the feasible region defined by the inequality (\ref{fs2}), i.e., 
\[\bar{F} = \{x\in \mathbb{R}^n_{+}\ |\ \sum_{i=1}^{n} \{\min_{y\in Y} a_i(y)\} x_i + \min_{y\in Y} a_0(y)\geq 0\}.\]
From above analysis, we conclude that $\bar{F} \subseteq F$.

The main difference between the original constraint (\ref{fs1}) and the restriction constraint (\ref{fs2}) is that
the minimization is independent on the decision variable $x$ in the later case.
In order to compute $\bar{F}$, we need to solve a series of problems as follows:
\begin{equation}
\label{gloa}
\min_{y} \ a_i(y)\quad \textrm{s.t.}\quad y\in Y
\end{equation}
for $i = 0,1,...,n$.
Based on $\bar{F}$, we can construct a linear program associated with one linear 
inequality constraint such that it has the same objective
function as LSIP and any feasible point of the constructed problem is feasible for LSIP. 
Such a problem is defined as
\begin{equation}
  \min_{x\in \mathbb{R}^n} \quad c^{\top} x \quad 
   \textrm{s.t.} \quad x \in \bar{F}.
   \tag*{(R-LSIP)}\label{rlsip}
\end{equation}
To characterize how well R-LSIP approximates LSIP, we can estimate the distance between
$g(x) =  \min_{y\in Y}  \{a(y)^{\top} x + a_0(y)\}$ and 
$\bar{g}(x) = \sum_{i=1}^{n} \{\min_{y\in Y} a_i(y)\} x_i +
\min_{y\in Y} a_0(y)$ which have been used to define the constraints of LSIP and R-LSIP.
Assume that each function $a_i(y)$ is Lipschitz continuous on $Y$, i.e., there exist some 
constant $L_i \geq 0$ such that $|a_i(y) - a_i(z)| \leq L_i |y - z|$ holds for any 
$y,z\in Y, i = 0,1,2,...,n$. By direct computation, we have 
\[|g(x) - \bar{g}(x)| \leq (\sum_{i=1}^{n} L_i x_i + L_0) (b-a).\]
It turns out that for any fixed $x$, the error between $g(x)$ and $\bar{g}(x)$ is bounded  
linearly with respect to $(b-a)$. 
Furthermore, if we assume that the decision variables are upper bounded 
(e.g., $0\leq x_i \leq U_i$ for some constants $U_i > 0$, $i = 0,1,2,...,n$), we have
\[|g(x) - \bar{g}(x)| \leq (\sum_{i=1}^{n} L_i U_i + L_0) (b-a).\]
This indicates that the error between $g(x)$ and $\bar{g}(x)$ goes to zeros uniformly
as $|b-a|$ tends to zero. By dividing the index set $Y = [a,b]$ into subintervals, 
one can construct a sequence of linear programs that approximate LSIP exhaustively 
as the size of the subdivision (formally defined in section {\bf 3.1}) tends to zero. 
Given a subdivision, constructing R-LSIP on each subinterval requires to solve \ref{rlsip} 
globally which will become computationally expensive due to the increasing number of 
subintervals and non-convexity of the coefficient functions in general. In fact, it is not 
necessary to solve \ref{rlsip} exactly. In the next section, we will discuss how to  
estimate a good lower bound of \ref{rlsip} and use it to construct the feasible approximation
problems for \ref{LSIP}.

\subsection{Construction of the inner approximation region}
In order to guarantee that the feasible region $\bar{F}$ derived from inequality (\ref{fs2}) is an
inner approximation of the feasible region of \ref{LSIP}, optimization problem (\ref{gloa}) 
needs to be solved globally. However, computing a lower bound for (\ref{gloa}) is enough to generate 
a restriction problem of \ref{LSIP}. In this section, we present two alternative approaches to 
approximate problem (\ref{gloa}).  
The idea of the first approach comes from the techniques of interval methods \cite{IT,IT1}.
Given an interval $Y = [a,b]$, the range of $a_i(y)$ on $Y$ is defined as $R(a_i,Y) = [R_i^l,R_i^u] = \{a_i(y)\ |\ y\in Y\}$. An interval function $A_i(Y) = [A_i^l, A_i^u]$ is called a inclusion function for $a_i(y)$ on $Y$ if $R(a_i,Y) \subseteq A_i(Y) $.
A natural inclusion function can be obtained by replacing the decision variable $y$ in $a_i(y)$ with the corresponding
interval and computing the resulting expression using the rules of interval arithmetic \cite{IT1}.
In some special cases, the natural inclusion function is tight (i.e., $R(a_i,Y) = A_i(Y)$). However, 
in more general cases, the natural interval function overestimates the original range of 
$a_i(y)$ on $Y$ which implies that $A_i^l < \min_{y\in Y} a_i(y)$.  In such cases, the tightness
of the inclusion can be measured by
\begin{equation}\label{haus}
  \max\{|R_i^l - A_i^l|,|R_i^u - A_i^u|\} \leq \gamma |b-a|^p\quad \textrm{and}\quad |A_i^l - A_i^u| \leq \delta |b-a|^p,
\end{equation}
where $p \geq 1$ is the convergence order, $\gamma \geq0$ and $\delta\geq 0$ are constants which 
depend on the expression of $a_i(y)$ and the interval $[a,b]$.
By replacing $\min_{y\in Y} a_i(y)$ in (\ref{fs2}) with $A_i^l$ for $i = 0,1,...,n$, we have a new 
linear inequality as follows
\begin{equation}\label{app1}
  \sum_{i = 1}^{n}A_i^l x_i + A_0^l \geq 0.
\end{equation}
It is obvious that any $x$ satisfying (\ref{app1}) is a feasible point for \ref{LSIP}.

The second approach to estimate the lower bound of problem (\ref{fs2}) is to construct a uniform 
lower bound function $\bar{a}_i(y)$  such that $\bar{a}_i(y) \leq a_i(y)$ holds for all 
$y\in Y$. In addition, we require that the optimal solution for 
\[\min_{y}\ \bar{a}_i(y)\quad \textrm{s.t.}\quad y\in Y\]
is easy to be identified. Here, we construct a concave lower bound function for $a_i(y)$ by adding a 
negative quadratic term to it, i.e., 
\[\bar{a}_i(y) = a_i(y) - \frac{\alpha_i}{2}(y-\frac{a+b}{2})^2,\]
where $\alpha \geq 0$ is a parameter. It follows that $\bar{a}_i(y) \leq a_i(y) \ \forall y\in Y$.
Furthermore, $\bar{a}_i(y)$ is twice continuously differentiable if and only if
$a_i(y)$ is twice continuously differentiable and the second derivative of $\bar{a}_i(y)$ is
$\bar{a}''_i(y) = a''_i(y) - \alpha_i$.
Thus $\bar{a}_i(y)$ is concave on $Y$ if the parameter $\alpha_i$ satisfies
$\alpha_i \geq \max_{y\in Y}\ a''_i(y)$. To sum up, we select the parameter $\alpha_i$ such that
\begin{equation}\label{alpha}
  \alpha_i \geq \max \{0, \max_{y\in Y} a''_i(y)\}.
\end{equation}
This guarantees that $\bar{a}_i(y)$ is a  lower bound concave function of $a_i(y)$ on the index set
$Y$. The computation of $\alpha_i$ in (\ref{alpha}) involves a global optimization. However, we can use 
any upper bound of the right hand side in (\ref{alpha}). Such an upper bound can be obtained by interval methods proposed above.
On the other hand, the distance between $\bar{a}_i(y)$ and $a_i(y)$ on $[a,b]$ is
\begin{equation*}\label{dist}
  \max_{y\in Y} |a_i(y) - \bar{a}_i(y)| = \frac{\alpha_i}{8}(b-a)^2.
\end{equation*}
Since $\bar{a}_y(y)$ is concave on $Y$, the minimizer of $\bar{a}_i(y)$ on $Y$ is attained on the boundary of $Y$ (see \cite{convex}), i.e., 
$\min_{y\in Y} \bar{a}_i(y) = \min\{\bar{a}_i(a),\bar{a}_i(b)\}$.
By replacing $\min_{y\in Y} a_i(y)$ in (\ref{fs2}) with $\min_{y\in Y} \bar{a}_i(y)$, we get the 
second type of restriction constraint as follows 
\begin{equation}\label{app2}
  \sum_{i = 1}^{n} \min\{\bar{a}_i(a),\bar{a}_i(b)\} x_i + \min\{\bar{a}_0(a),\bar{a}_0(b)\} \geq 0.
\end{equation}

The two approaches are distinct in the sense that the interval method requires mild assumptions on
the coefficient function while the concave-function based method admits better approximation rate.


\section{Numerical method}
Based on the restriction approaches developed in the previous section, we are able to construct a 
sequence of approximations for \ref{LSIP} by dividing the original index set into subsets 
successively and constructing linear optimization problems associated with restricted constraints on 
the subsets. 
\begin{definition}
We call $T = \{\tau_0,...,\tau_N\}$ a subdivision of the interval $[a,b]$ if
\[a=\tau_0 \leq \tau_1\leq ... \leq \tau_N = b.\]
\end{definition}
Let $Y_k = [\tau_{k-1},\tau_{k}]$ for $k = 1,2,...,N$, the length of $Y_k$ is defined by $|Y_k| = |\tau_k - \tau_{k-1}|$ and the
length of the subdivision $T$ is defined by $|T| = \max_{1\leq k \leq N}|Y_k|$.
It follows that $Y = \cup_{i = 1}^N Y_k$. 

The intuition behind the approximation of \ref{LSIP} through subdivision comes from an observation 
that the original semi-infinite constraints in \ref{LSIP}
\[a(y)^{\top}x + a_0(y) \geq 0, \ \forall y\in Y\]
can be reformulated equivalently as finitely many semi-infinite constraints
\[a(y)^{\top}x + a_0(y) \geq 0, \ \forall y\in Y_k, k = 1,2,...,N.\]
Given a subdivision, we can construct the approximate constraint on each subinterval and combine
them together to formulate the inner-approximation of the original feasible region. The corresponding
optimization problem provide a restriction of \ref{LSIP}. The solution of the approximate problem
approach to the optimal solution of \ref{LSIP} as the size of the subdivision tends to zero.

The two different approaches (e.g., interval method and Concavification method)
were introduced in section 2 to construct the approximate region  that lies inside of the original 
feasible region. This induces two different types of approximation problems when applied to a 
particular subdivision. We only describe main results for the first type (e.g., interval method) 
and focus on the convergence and algorithm for the second one.

We introduce the Slater condition and a lemma derived from it which will be used in the
following part. We say Slater condition holds for \ref{LSIP} if there exists a point $\bar{x}\in \mathbb{R}^n_{+}$ such that
\[a(y)^{\top}\bar{x}+a_0(y) > 0, \ \forall y\in Y.\]
Let $F^o = \{x\in F\ |\ a(y)^{\top}\bar{x}+a_0(y) > 0, \ \forall y\in Y\}$ be the set of all the Slater points in $F$.
It is shown that the feasible region $F$ is exactly the closure of $F^o$
under the Slater condition \cite{mma}. We present this result as a lemma and give a direct proof in 
the appendix.
\begin{lemma}\label{lm32}
Assume that the Slater condition holds for \ref{LSIP} and the index set $Y$ is compact, 
then we have
\[ F = cl(F^o),\]
where $cl(F^o)$ represents the closure of the set $F^o$.
\end{lemma}

\subsection{Restriction based on interval method}
Let $A_i(Y_k) = [A_{i,k}^l,A_{i,k}^u]$ be the inclusion function of $a_i(y)$ on $Y_k$. By estimating
the lower bound for $\min_{y\in Y_k} a_i(y)$ via interval method, we can construct the following
linear constraints 
\[  \sum_{i=1}^n A_{i,k}^l x_i + A_{0,k}^l \geq 0, k = 1,2,...,N,\]
corresponding to the original constraints 
$a(y)^{\top}x + a_0(y) \geq 0, \ \forall y\in Y_k, k = 1,2,...,N$.
For simplicity, we reformulate the inequalities as 
\begin{equation}
A_T^{\top}x+b_T \geq 0,
\label{31}
\end{equation}
where $A_T(i,k) = A_{i,k}^l$ and $b_T(k) = A_{0,k}^l$ for $i = 1,2,...,n,  k = 1,2,...,N$.
The approximation problem for \ref{LSIP} in such case is formulated as 
\begin{align}
\begin{split}
\min_{x\in \mathbb{R}^n_{+}}\ c^{\top}x \quad \textrm{s.t.}\quad A_T^{\top}x+b_T \geq 0.
\end{split}\tag*{R1-LSIP(T)}\label{r1lsip}
\end{align}
Following the analysis in section 2, we know that
$\{x\in \mathbb{R}^n_{+}\ |\ A_T^{\top}x+b_T \geq 0\} \subseteq F$.
Therefore, any feasible point of \ref{r1lsip} is feasible for \ref{LSIP} provided that the feasible 
region of \ref{r1lsip} is non-empty. By solving \ref{r1lsip}, we can obtain a feasible approximate solution for \ref{LSIP} and the corresponding optimal value
of \ref{r1lsip} provides an upper bound for the optimal value of \ref{LSIP}.

Let $F(T) = \{x\in \mathbb{R}^n_{+}\ |\ A_T^{\top}x + b_T \geq 0\}$ be 
the feasible region of \ref{r1lsip}. We say that $F(T)$ is consistent if
$F(T) \neq \emptyset$. In this case, the corresponding problem \ref{r1lsip} is called consistent.
The following lemma shows that the approximate problem \ref{r1lsip} is consistent for all $|T|$ small enough if Slater condition holds for \ref{LSIP}.

\begin{lemma}\label{lm33}
  Assume that the Slater condition holds for \ref{LSIP} and the coefficient functions 
  $a_i(y)$, $i = 0,1,...n$, are
  Lipschitz continuous on $Y$, then $F(T)$ is nonempty for all $|T|$ small enough.
\end{lemma}

In following theorem, we show that any accumulation point of the solutions of the approximate 
problems \ref{r1lsip} is a solution to \ref{LSIP} if the size of the subdivision tends to zero.
\begin{theorem}\label{t34}
Assume the Slater condition holds for \ref{LSIP} and the level set
$L(\bar{x}) = \{x\in F\ |\ c^{\top}x\leq c^{\top}\bar{x}\}$ is bounded ($\bar{x}$ is a Slater
point). Let $\{T_k\}$ be a sequence of subdivisions of $Y$ such that
$T_0$ is consistent and $\lim_{k\to \infty} |T_k| = 0$ with $T_k\subseteq T_{k+1}$. 
Let $x_k^*$ be a solution of R1-LSIP($T_k$). Then any accumulation point of the sequence $\{x_k^* \}$ is an optimal solution to \ref{LSIP}.
\end{theorem}


\subsection{Restriction based on concavification}
Given a subdivision $T = \{\tau_0,...,\tau_N\}$ and $Y_k = [\tau_{k-1},\tau_k], k=1,2,...N$, 
by applying concavification method in section 2 to each of the finitely many semi-infinite constraints
\[a(y)^{\top}x + a_0(y)\geq 0, \ \forall y\in Y_k, k = 1,2,...,N,\]
we can construct the linear constraints as follows
\[\sum_{i=1}^n \min\{\bar{a}_i(\tau_{k-1}),\bar{a}_i(\tau_k)\} x_i + \min\{\bar{a}_0(\tau_{k-1}),\bar{a}_0(\tau_k)\} \geq 0, \ k = 1,2,...,N,\]
where $\bar{a}_i(\cdot)$ is the concavification function defined on $Y_k$ when we calculate
$\bar{a}_i(\tau_{k-1})$ or $\bar{a}_i(\tau_k)$ (i.e., 
$\bar{a}_i(y) = a_i(y) - \frac{\alpha_{i,k}}{2} (y- \frac{\tau_k - \tau_{k-1}}{2})^2$).
We rewrite the above inequalities as
\[\bar{A}^{\top}_{T} x + \bar{b}_T \geq 0,\]
where $\bar{A}_T(i,k) = \min\{\bar{a}_i(\tau_{k-1}),\bar{a}_i(\tau_k)\}$ and 
$\bar{b}_T(k) = \min\{\bar{a}_0(\tau_{k-1}),\bar{a}_0(\tau_k)\}$. The corresponding approximate 
problem for \ref{LSIP} is defined by
\begin{align}
\begin{split}
\min_{x\in \mathbb{R}^n_{+}}\ c^{\top}x \quad \textrm{s.t.}\quad \bar{A}_T^{\top}x+\bar{b}_T \geq 0.
\end{split}\tag*{R2-LSIP(T)}\label{r2lsip}
\end{align}
Let $\bar{F}(T) = \{x\in \mathbb{R}^n_{+} \ |\ \bar{A}^{\top}_{T} x + \bar{b}_T \geq 0\}$ be the 
feasible set of the problem \ref{r2lsip}, we can conclude that $\bar{F}(T)\subseteq F$.

The approximate problem \ref{r2lsip} is similar to \ref{r1lsip} in the sense that both problems
induce restrictions of \ref{LSIP}. Therefore, any feasible solution of \ref{r2lsip} is feasible 
for \ref{LSIP} and the corresponding optimal value provide an upper bound for the optimal value
of the problem \ref{LSIP}.

The following lemma shows that if the Slater condition holds for \ref{LSIP}, 
\ref{r2lsip} is consistent for all $|T|$ small enough (e.g., $\bar{F}(T) \neq \emptyset$). 
Proof can be found in appendix.
\begin{lemma}\label{l35}
  Assume the Slater condition holds for \ref{LSIP} and $a_i(y)$, $i = 1,2,...,n$, are twice continuously differentiable. Then \ref{r2lsip} is consistent for all $|T|$ small enough.
\end{lemma}

In order to find a good approximate solution for \ref{LSIP}, \ref{r2lsip} need to be solved 
iteratively during which the subdivision will be refined. We present a particular strategy of 
the refinement here such that the approximate regions of \ref{r2lsip} are monotonically 
enlarging from the inside of the feasible region $F$. Consequently, the corresponding optimal
values of the approximation problems are monotonically decreasing and converge to the optimal
value of the original linear semi-infinite problem. Note that such a refinement procedure can
not guarantee the monotonic property when applied to solve \ref{r1lsip}.  

Let $T = \{\tau_k\ |\ k = 0,1,...,N\}$ be a subdivision of the $Y$. Assume 
$Y_k = [\tau_{k-1},\tau_k]$ is the subinterval to be refined.
Denote by $\tau_{k,1}$ and $\tau_{k,2}$ the trisection points of $Y_k$:
\[\tau_{k,1} = \tau_{k-1} + \frac{1}{3}(\tau_k - \tau_{k-1}), \ 
\tau_{k,2} = \tau_{k-1} + \frac{2}{3}(\tau_k -
\tau_{k-1}).\]
The constraint in \ref{r2lsip} on the subset $Y_k$ is
\begin{equation}\label{3.3}
  \sum_{i=1}^n [\min [\bar{a}_i(\tau_{k-1}), \bar{a}_i(\tau_k)]x_i + \min [\bar{a}_0(\tau_{k-1}), \bar{a}_0(\tau_k)] \geq 0,
\end{equation}
where $\bar{a}_i(y) = a_i(y) - \frac{\alpha_{i,k}}{2}(y - \frac{\tau_{k-1} + \tau_k}{2})^2$ and 
parameter $\alpha_{i,k}$ is calculated in the manner of (\ref{alpha}).
The lower bounding functions on each subset after refinement are defined by
\begin{align*}
  \bar{a}_i^1(y) &= a_i(y) - \frac{\alpha^1_{i,k}}{2}(y-\frac{\tau_{k-1}+\tau_{k,1}}{2})^2,\ y\in Y_{k,1}=[\tau_{k-1},\tau_{k,1}], \\
  \bar{a}_i^2(y) &= a_i(y) - \frac{\alpha^2_{i,k}}{2}(y-\frac{\tau_{k,1}+\tau_{k,2}}{2})^2,\ y\in Y_{k,2}=[\tau_{k,1},\tau_{k,2}], \\
  \bar{a}_i^3(y) &= a_i(y) - \frac{\alpha^3_{i,k}}{2}(y-\frac{\tau_{k,2}+\tau_{k}}{2})^2,\ y\in Y_{k,3}=[\tau_{k,2},\tau_{k}],
\end{align*}
where $\alpha_{i,k}^j, j=1,2,3$ are selected such that
$\alpha_{i,k}^j \geq \max \{0, \max_{y\in Y_{k,j}} \nabla^2 a_i(y)\}$ and $\alpha_{i,k}^j \leq \alpha_{i,k}$ for $j = 1,2,3$.
The refined approximate region $\bar{F}(T\cup\{\tau_{k,1},\tau_{k,2}\})$ is obtained by replacing
the constraint (\ref{3.3}) in $\bar{F}(T)$ with
\begin{align*}
  & \sum_{i=1}^n [\min [\bar{a}_i^1(\tau_{k-1}), \bar{a}_i^1(\tau_{k,1})]x_i + \min [\bar{a}_0^1(\tau_{k-1}), \bar{a}_0^1(\tau_{k,1})] \geq 0, \\
  & \sum_{i=1}^n [\min [\bar{a}_i^2(\tau_{k,1}), \bar{a}_i^2(\tau_{k,2})]x_i + \min [\bar{a}_0^2(\tau_{k,1}), \bar{a}_0^2(\tau_{k,2})] \geq 0, \\
  & \sum_{i=1}^n [\min [\bar{a}_i^3(\tau_{k,2}), \bar{a}_i^3(\tau_k)]x_i + \min [\bar{a}_0^3(\tau_{k,2}), \bar{a}_0^3(\tau_k)] \geq 0.
\end{align*}
\begin{lemma}\label{3.6}
  Let $T$ be a consistent subdivision of $Y$. Assume that $\bar{F}(T\cup\{\tau_{k,1},\tau_{k,2}\})$ 
  is obtained by the trisection refinement procedure above, then we have
  \[\bar{F}(T) \subseteq \bar{F}(T\cup\{\tau_{k,1},\tau_{k,2}\}) \subseteq F.\]
\end{lemma}
\begin{proof}
  Since $x\in \mathbb{R}^n_{+}$, it suffices to prove that for $i = 0,1,2,...,n,$
  \[\min[\bar{a}_i^1(\tau_{k-1}), \bar{a}_i^1(\tau_{k,1}),\bar{a}_i^2(\tau_{k,1}), \bar{a}_i^2(\tau_{k,2}),\bar{a}_i^3(\tau_{k,2}), \bar{a}_i^3(\tau_k)] \geq \min [\bar{a}_i(\tau_{k-1}), \bar{a}_i(\tau_k)].\]
  By direct computation, we know $\bar{a}_i^1(\tau_{k-1}) \geq \bar{a}_i(\tau_{k-1})$ and
  $\bar{a}_i^3(\tau_k) \geq \bar{a}_i(\tau_k)$. Since $\bar{a}_i(y)$ is concave
  on $Y_k = [\tau_{k-1},\tau_k]$, we have $\bar{a}_i(\tau_{k,j}) \geq \min [\bar{a}_i(\tau_{k-1}), \bar{a}_i(\tau_k)]$ for $j = 1,2$. 
  In addition, direct calculation implies 
  \begin{equation*}
    \min[\bar{a}_i^1(\tau_{k,1}),\bar{a}_i^2(\tau_{k,1})] \geq \bar{a}_i(\tau_{k,1}), \quad
    \min[\bar{a}_i^2(\tau_{k,2}),\bar{a}_i^3(\tau_{k,2})] \geq \bar{a}_i(\tau_{k,2}).
  \end{equation*}
  The last two statements indicate that
$    \min[\bar{a}_i^1(\tau_{k,1}),\bar{a}_i^2(\tau_{k,1})] \geq \min [\bar{a}_i(\tau_{k-1}), \bar{a}_i(\tau_k)]$ and
$ \min[\bar{a}_i^2(\tau_{k,2}),\bar{a}_i^3(\tau_{k,2})] \geq \min [\bar{a}_i(\tau_{k-1}), \bar{a}_i(\tau_k)]$. 

This proves our statement.
\qquad\end{proof}

We present in the following theorem the general convergence results for approximating \ref{LSIP} via a 
sequence of restriction problems.
\begin{theorem}\label{3.7}
  Assume that the assumptions in Theorem \ref{t34} hold.
  Let $\{T_k\}$ be a sequence of subdivisions of the index set $Y$, which is obtained by trisection refinement recursively, such that $T_0$ is consistent and $\lim_{k\to\infty} |T_k| = 0$.
  Denote by $x^*_k$ the optimal solution to R2-LSIP($T_k$). Then we have:\\
 (1) $x^*_k$ is feasible for \ref{LSIP} and any accumulation point of the sequence $\{x^*_k\}$
    is a feasible solution to \ref{LSIP}.\\
(2) $\{f(x^*_k):\ f(x^*_k) = c^{\top}x^*_k\}$ is a decreasing sequence and
    $v^* = \lim_{k\to\infty} f(x^*_k)$ is an optimal value to \ref{LSIP}.
\end{theorem}
\begin{proof}
The proof of the first statement is similar to the proof in Theorem \ref{t34}.

From Lemma \ref{3.6}, we know that $\bar{F}(T_{k-1}) \subseteq \bar{F}(T_{k})$ holds 
for $k\in \mathbb{N}$ which implies that the sequence $\{f(x^*_k)\}$ is decreasing.
Since the level set $L(\bar{x})$ is bounded, the sequence $\{f(x^*_k)\}$ is bounded. Therefore,
the limit of the sequence exists which is denoted by $v^*$. From (1), we know that
$v^*$ is an optimal value to \ref{LSIP}.

  This completes our proof.
\qquad\end{proof}


\subsection{Adaptive refinement algorithm}
In this section, we present a specific algorithm to solve \ref{LSIP}. The algorithm is based on 
solving the approximate linear problems \ref{r2lsip}  (or \ref{r1lsip}) for a given subdivision 
$T$ and then refine the subdivision to improve the solution. The key idea of the algorithm is to
select the candidate subsets in $T$ to be refined in an adaptive manner rather than making the
refinement exhaustively. 

We introduce the optimality condition for \ref{LSIP} as follows before presenting the details of 
the algorithm.
Given a point $x\in F$, let $A(x) = \{y\in Y\ |\ a(y)^{\top}x + a_0(y) = 0\}$ be the active index set
for \ref{LSIP} at $x$. If some constraint qualification (e.g., Slater condition) holds for \ref{LSIP},
a feasible point $x^* \in F$ is an optimal solution if and only if $x^*$ satisfies the KKT systems
(\cite{mma}), i.e.,
\[c - \sum_{y\in A(x^*)}\lambda_y a(y) = 0\]
for some $\lambda_y \geq 0$, $y\in A(x^*)$.

\begin{definition}
  We say that $x^*\in F$ is an $(\epsilon,\delta)$ optimal solution to \ref{LSIP} if
  there exist some indices $y\in Y$ as well as $\lambda_y \geq 0$ such that
  \[||c - \sum_{y\in A(x^*,\delta)}\lambda_y a(y)||\leq \epsilon,\]
  where $A(x^*,\delta) = \{y\in Y\ |\ 0\leq a(y)^{\top}x^* + a_0(y) \leq \delta\}$.
\end{definition}

\begin{table*}
\indent {\bf ------------------------------------------------------------------------------------------------}\\
\indent\textbf{Algorithm 1} (Adaptive Refinement Algorithm for LSIP) \\
\indent {\bf ------------------------------------------------------------------------------------------------}
\begin{itemize}
  \item[\textbf{S1.}] Find an initial subdivision $T_0$ such that R2-LSIP($T_0$) is consistent. 
  Choose an initial point $x_0$ and tolerances $\epsilon$ and $\delta$. Set $k = 0$.
  \item[\textbf{S2.}] Solve R2-LSIP($T_k$) to obtain a solution $x^*_k$
   and the active index set $A(x^*_k)$.
  \item[\textbf{S3.}] Terminate if $x^*_k$ is an $(\epsilon,\delta)$ optimal solution to 
  \ref{LSIP}.
  Otherwise update $T_{k+1}$ and $F(T_{k+1})$ by trisection refinement procedure for subintervals 
  in $T_k$ that correspond to $A(x^*_k)$.
  \item[\textbf{S4.}] Let $k = k+1$ and go to  step 2.
\end{itemize}
\indent {\bf ------------------------------------------------------------------------------------------------}
\end{table*}

To obtain a consistent subdivision in the firs step of Algorithm 1, we apply the adaptive 
refinement algorithm to the following problem
\begin{align}\label{initial}
\begin{split}
  \min_{(x,z)\in \mathbb{R}^n_{+}\times \mathbb{R}}\ z \quad \textrm{s.t.}\quad a(y)^{\top}x + a_0(y) \geq z\ \forall y\in Y
\end{split}\tag*{LSIP$_0$}
    \end{align}
until a feasible solution $(x_0,z_0)$, with $z_0 \geq 0$, of the problem LSIP$_0$($T_0$) is found 
for some subdivision $T_0$. The current subdivision $T_0$ is consistent and chosen as the initial 
subdivision of Algorithm 1. In addition, $x_0$ is feasible for the original problem and selected 
as the initial point for the algorithm.

The refinement procedure in the third step of the algorithm is taken as follows.
In the $k$th iteration, each $[\tau^k_{i-1},\tau^k_i]$ is divided into three equal length subsets 
for $i\in A(x^*_k)$. New constraints are constructed on the subsets and used to update the 
constraint corresponding to $[\tau^k_{i-1},\tau^k_i]$ for each index $i\in A(x^*_k)$.
Then we have $\bar{F}(T_{k+1})$ and the associated approximation problem R2-LSIP($T_{k+1}$).


\begin{theorem}[Convergence of Algorithm 1]
\label{conapp2}
Assume the Slater condition holds for \ref{LSIP} and the coefficient functions 
$a_i(y)$, $i = 0,1,...,n$, are twice continuously differentiable. Then Algorithm 1
terminates in finitely many iterations for any positive tolerances $\epsilon$ and $\delta$.
\end{theorem}
\begin{proof}
  Let $x^*_k$ be a solution to the approximate subproblem R2-LSIP($T_k$) with 
  $T_k = \{\tau^k_j\ |\ j = 0,1,...,N_k\}$, there exists some
  $\lambda_j^k \geq 0$ for $j\in A(x^*_k)$ such that
  \begin{equation}\label{kkt}
    c - \sum_{j\in A(x^*_k)} \lambda_j^k \min[\bar{a}(\tau^k_{j-1}),\bar{a}(\tau_j^k)] = 0,
  \end{equation}
  where $A(x^*_k) = \{j\ |\ \min[\bar{a}(\tau^k_{j-1}),\bar{a}(\tau_j^k)]x^*_k +
   \min[\bar{a}_0(\tau^k_{j-1}),\bar{a}_0(\tau_j^k)]= 0\}$ is the active index set for R2-LSIP($T_k$)
  at $x^*_k$ and $\min[\bar{a}(\tau^k_{j-1}),\bar{a}(\tau_j^k)]$ represents a vector in $\mathbb{R}^n$ 
  such that the $i$th element is defined by $\min[\bar{a}_i(\tau^k_{j-1}),\bar{a}_i(\tau_j^k)]$.
  Since $\bar{a}_i(y) = a_i(y)-\frac{\alpha_{i,k}}{2}(y-\frac{\tau_{j-1}^k + \tau_j^k}{2})^2$ for
  $y\in [\tau^k_{j-1},\tau^k_j]$, we have
   \[\min[\bar{a}(\tau^k_{j-1}),\bar{a}(\tau_j^k)] = \min[a(\tau^k_{j-1}),a(\tau_j^k)] - \frac{1}{8}(\tau^k_j -
   \tau^k_{j-1})^2 \alpha^k,\]
   where $\alpha^k = (\alpha^k_{1,j},\alpha^k_{2,j},...,\alpha^k_{n,j})^{\top}$ is the parameter vector on the subset $[\tau^k_{j-1},\tau^k_j]$ with all elements are uniformly bounded. On the other hand, since $a_i(y)$ is twice continuously differentiable, there exists $\bar{\tau}_{j-1}^k$
   such that 
   \[a_i(\tau^k_j) = a_i(\tau^k_{j-1}) + a_i^{'}(\bar{\tau}^k_{j-1})(\tau^k_j - \tau^k_{j-1}), 1\leq i\leq n\]
   which implies that $\min[a(\tau^k_{j-1}),a(\tau_j^k)] = a(\tau_{j-1}^k) + 
   (\tau^k_j - \tau^k_{j-1}) \beta^k$ where $\beta^k \in \mathbb{R}^n$ is a constant vector
   (e.g., $\beta^k_i = a_i^{'}(\bar{\tau}^k_{j-1})$ if $\min[a(\tau^k_{j-1}),a(\tau_j^k)] 
   = a(\tau_j^k)$ and $\beta^k_i = 0$ otherwise).
   It follows that 
   \begin{equation}
   \min[\bar{a}(\tau^k_{j-1}),\bar{a}(\tau_j^k)] = a(\tau^k_{j-1}) + 
   (\tau^k_j - \tau^k_{j-1}) \beta^k - \frac{1}{8}(\tau^k_j - \tau^k_{j-1})^2 \alpha^k.
   \label{minabar}
   \end{equation}
Substitute $\min[\bar{a}(\tau^k_{j-1}),\bar{a}(\tau_j^k)]$ in (\ref{minabar}) into (\ref{kkt}) and
$A(x^*_k)$, we can claim that it suffices to prove the lengths of all the subsets 
$[\tau^k_{j-1},\tau^k_j]$ for $j\in A(x^*_k)$ converge to zeros as the iteration $k$ tends to infinity.
From the algorithm, we know that in each iteration at least one subset $[\tau^k_{j-1},\tau^k_j]$ is
divided into three equal subintervals where the length of each subinterval is bounded above by 
$\frac{1}{3}(\tau^k_j - \tau^k_{j-1}) \leq \frac{1}{3}(b-a)$. 
For each integer $p\in \mathbb{N}$, at least one interval with its length bounded by
$\leq \frac{1}{3^p}(b-a)$ is generated. Furthermore, all the subintervals 
$[\tau^k_{j-1},\tau^k_j]$, $j\in A(x^*_k)$ are different for all $k\in \mathbb{N}$. Since for each 
$p\in \mathbb{N}$, only finitely subintervals with length greater than $\frac{1}{3^p}(b-a)$ exists. 
This implies that the lengths of the subsets
$[\tau^k_{j-1},\tau^k_j]$ for $j\in A(x^*_k), k\in \mathbb{N}$ must tend to zero.
\qquad\end{proof}

We can conclude from Theorem \ref{conapp2} that if the tolerances $\epsilon$ and $\delta$ are 
decreasing to zero then any accumulation point of the sequence generated by Algorithm 1 is a solution 
to the original linear semi-infinite programming. 

\begin{corollary}
Let the assumptions in Theorem \ref{conapp2} be satisfied and the tolerances $(\epsilon_k,\delta_k)$ 
are chosen such that $(\epsilon_k,\delta_k)\searrow (0,0)$. If $x^*_k$ is an $(\epsilon_k,\delta_k)$ KKT point for \ref{LSIP} generated by Algorithm 1, then any accumulation point $x^*$ of the sequence $\{x^*_k\}$ is a solution to \ref{LSIP}.\label{corollary310}
\end{corollary}

It follows from Corollary \ref{corollary310} the sequence $\{c^{\top}x^*_k\}$ is monotonically 
decreasing to the optimal value of \ref{LSIP}  as $k$ tends to infinity. In the implement of our
algorithm, the termination criterion is set as
\[|c^{\top}x^*_k - c^{\top}x^*_{k-1}| \leq \epsilon.\]


The convergence of Algorithm 1 is also applicable to the case that the approximate problem 
R1-LSIP($T_k$) is used in the second step. The proof is similar to that in theorem \ref{conapp2} as we 
explained in appendix. However, we can not guarantee the sequence $\{c^{\top}x^*_k\}$  is 
monotonically decreasing.


\subsection{Remarks}
The proposed algorithm can be applied to solve linear semi-infinite optimization problem with finitely 
many semi-infinite constraints and some extra linear constraints, i.e., 
\[\min_{x\in X} c^{\top}x \quad \textrm{s.t.} \quad a^j(y)^{\top}x + a^j_0(y) \geq 0, 
\ \forall y\in Y, j = 1,2,...,m,\]
where $X = \{x\in \mathbb{R}^n\ |\ Dx\geq d\}$ and $a^j(\cdot): \mathbb{R}
\mapsto \mathbb{R}^n$.
 In such a case, we split each decision variable  $x_i$ into two non-negative variable $y_i \geq 0$ 
and $z_i \geq 0$ such that $x_i = y_i - z_i$, and then substitute $x_i$ into the above problem.
Then the problem is reformulated as a linear semi-infinite programming problem with non-negative 
decision variables in which the Algorithm 1 can be applied to solve it. 
Such a technique is applied in the numerical experiments.

In the case that $X = [X_l, X_u]$ is a box in $\mathbb{R}^n$, we can set a new variable transformation
as $x = z + X_l$ in which $z\geq 0$. 
The advantage to reformulate the original problem in such a translation is that the
dimension of the new variables is the same as that of the original decision variables.


\section{Numerical experiments}
We present the numerical experiments for a couple of optimizaiton problems selected from the 
literature. The algorithm is implemented
in $Matlab$ 8.1 and the subproblem is solved by using $linprog$ of $Optimization$ $Toolbox$ 6.3 with 
default tolerance and active set algorithm. All the following experiments were run on 3.2 GHz Intel(R) 
Core(TM) processor.

The computation of the bounds for the coefficient functions and the parameter $\alpha$ in the
second approach are obtained directly if the closed form bound exists. Otherwise, we use $Matlab$ 
toolbox $Intlab$ 6.0 \cite{intlab} to obtain the corresponding
bounding values. 
The problems in the literature are listed as follows. \\
\textbf{Problem\ 1.}
\begin{align*}
 \label{eg1}
   \min             \quad & \sum_{i=1}^n i^{-1}x_i  \\
    \textrm{s.t.}   \quad & \sum_{i=1}^n y^{i-1}x_i\geq tan(y), \  \forall y\in [0,1].
\end{align*}
This problem is taken from \cite{mma} and also tested in \cite{BB} for $n=8$. For $1\leq n\leq 7$, the problem
 has unique optimal solution and the Strong Slater condition holds. The problem for $n=8$ is hard to solve
 and thus a good test of the performance for our algorithm. \\
\textbf{Problem\ 2.} This problem has same formulation as Problem 1 with $n = 9$ which is also
tested in \cite{BB}. \\
\textbf{Problem\ 3.}
\begin{align*}
   \min             \quad & \sum_{i=1}^8 i^{-1}x_i  \\
   \textrm{s.t.}    \quad & \sum_{i=1}^8 y^{i-1}x_i\geq \frac{1}{2-y}, \  \forall y\in [0,1].
\end{align*}
This problem is taken from \cite{KS} and also tested in \cite{BB}. \\
\textbf{Problem\ 4.}
\begin{align*}
   \min             \quad & \sum_{i=1}^7 i^{-1}x_i  \\
   \textrm{s.t.}    \quad & \sum_{i=1}^7 y^{i-1}x_i\geq -\sum_{i=0}^4 y^{2i}, \  \forall y\in [0,1].
\end{align*}
\textbf{Problem\ 5.}
\begin{align*}
   \min             \quad & \sum_{i=1}^9 i^{-1}x_i  \\
   \textrm{s.t.}    \quad & \sum_{i=1}^9 y^{i-1}x_i\geq \frac{1}{1+y^2}, \  \forall y\in [0,1].
\end{align*}
Problem 4 and 5 are taken from \cite{TE} and also tested in \cite{BB}.

 The following problems, as noted in \cite{BB}, arise in the design of finite impulse response(FIR) filters which
 are more computationally demanding than the previous ones (see, e.g., \cite{rg,BB}).

\textbf{Problem\ 6.}
\begin{align*}
   \min             \quad & -\sum_{i=1}^{10} r_{2i-1}x_i  \\
   \textrm{s.t.}    \quad & 2\sum_{i=1}^{10} \textrm{cos}((2i-1)2\pi y)x_i \geq -1, \  \forall y\in [0,0.5],
\end{align*}
where $r_i = 0.95^i$.

\textbf{Problem\ 7.} This problem is formulated as Problem 6 where $r_i = 2\rho
\textrm{cos}(\theta)r_{i-1} - \rho^2
r_{i-2}$ with $\rho  = 0.975$, $\theta = \pi/3$, $r_0 = 1$, $r_1 = 2\rho \textrm{cos}(\theta)/(1+\rho^2)$.

\textbf{Problem\ 8.} This problem is also formulated as Problem 6 where $r_i = \frac{\textrm{sin}(2\pi f_s i)}{2\pi f_s i}$ with $f_s = 0.225$.

The numerical results are summarized in Table \ref{table1} where CPU Time is the time cost when 
the algorithm terminates, Objective Value represents the objective function value at the iteration
point when the algorithm terminates, No of Iteration is the number of 
iterations when the algorithm terminates for each particular problem and Violation measures the 
feasibility of the solution $x^*$ obtained by the algorithm which is defined by $\min_{y\in \bar{Y}} 
g(x^*,y)$ with $\bar{Y} = a:10^{-6}:b$.
We also list the numerical results for these problems by MATLAB toolbox $fseminf$ as a reference.
We can see that the algorithm proposed in this paper generates the feasible solutions  
for all the problems tested. This is coincide with the theoretical results. 
Furthermore, Algorithm 1 works well for the 
computational demanding problems 6-8. The solver $fseminf$ is faster than our method, however the 
feasibility is not guaranteed for this kind of method.

\renewcommand{\arraystretch}{1.3}
\begin{table}
\tbl{Summary of numerical results for the proposed algorithm in this paper}
{\begin{tabular}{lcccccc} \toprule
  & Algorithm    &   CPU Time(sec)  &   Objective Value  &   No. of Iterations & Violation \\ \midrule
\multirow{3}*{\bf Problem 1.}
           &  Approach 1    & 1.8382        & 0.6174   &      169        &  2.3558e-04  \\
           &  Approach 2    & 1.8910        & 0.6174   &      172        &  1.5835e-04  \\
           &  fseminf       & 0.2109        & 0.6163   &      33         &  -1.2710e-04  \\ \hline
\multirow{3}*{\bf Problem 2.}
           &  Approach 1    & 5.2691            & 0.6163   &      273       &  4.1441e-04 \\
           &  Approach 2    & 4.0928            & 0.6166   &      266       &  1.8372e-04 \\
           &  fseminf       & 0.3188            & 0.6157   &      46        &  -7.6194e-04 \\ \hline
\multirow{3}*{\bf Problem 3.}
           &  Approach 1    & 0.1646            & 0.6988   &       12       &  2.7969e-03 \\
           &  Approach 2    & 0.1538            & 0.6988   &       13       &  2.8014e-03 \\
           &  fseminf       & 0.2387            & 0.6932   &       35       &  -5.8802e-07 \\ \hline
\multirow{3}*{\bf Problem 4.}
           &  Approach 1    & 4.1606            & -1.7841  &      354      &  1.9689e-05 \\
           &  Approach 2    & 4.1928            & -1.7841  &      356      &  1.9646e-05 \\
           &  fseminf       & 0.4794            & -1.7869  &      70       &  -3.4649e-09  \\ \hline
\multirow{3}*{\bf Problem 5.}
           &  Approach 1    & 4.2124            & 0.7861   &      300  &  1.9829e-05 \\
           &  Approach 2    & 4.7892            & 0.7861   &      302  &  1.9243e-05 \\
           &  fseminf       & 0.3642            & 0.7855   &      32   &  -8.5507e-07   \\ \hline\hline
\multirow{3}*{\bf Problem 6.}
           &  Approach 1   & 1.7290            &  -0.4832  &       137      &    5.0697e-06 \\
           &  Approach 2   & 1.5302            &  -0.4832  &       132      &    5.0914e-06 \\
           &  fseminf      & 1.1476            &  -0.4754  &       86       &    -1.2219e-04 \\ \hline
\multirow{3}*{\bf Problem 7.}
           &  Approach 1   & 2.5183            &  -0.4889  &       170      &    2.8510e-04 \\
           &  Approach 2   & 3.2521            &  -0.4890  &       219      &    2.8861e-04 \\
           &  fseminf      & 1.0480            &  -0.4883  &       86       &    -1.5211e-03 \\ \hline
\multirow{3}*{\bf Problem 8.}
           &  Approach 1   & 4.4262            &  -0.4972  &       252      &    4.5808e-05 \\
           &  Approach 2   & 4.0216            &  -0.4972  &       252      &    5.0055e-05 \\
           &  fseminf      & 0.4324            &  -0.4973  &       45       &    -4.3322e-07 \\ \hline\hline
\end{tabular}}
\tabnote{\textsuperscript{*}
Approach 1 represents {\bf Algorithm 1} with R1-LSIP and Apporach 2 represents {\bf Algorithm 1} with
R2-LSIP.}
\label{table1}
\end{table}

\section{Conclusion}
A new numerical method for solving linear semi-infinite programming problems
is proposed which guarantees that each iteration point is feasible for the original problem. 
The approach is based on a two-stage restriction of the original semi-infinite constraint. 
The first stage restriction allows us to consider semi-infinite constraint independently to the 
decision variables on the subsets of the index set. In the second stage, 
the lower bounds for the optimal values of the optimization problems associated with coefficient 
functions are estimated using two different approaches. The approximation error goes to zero
as the size of the subdivisions tends to zero.

The approximate problems with finitely many linear constraints is constructed such that the
corresponding feasible regions are included in the feasible region of \ref{LSIP}. It follows that 
any feasible solution of the approximate problem is feasible for \ref{LSIP} and the corresponding 
objective function value provide an upper bound for the optimal value of \ref{LSIP}. It is proved 
that the solutions of the approximate problems converge to that of the original problem. 
Also, the sequence of optimal values of the approximate problems converge to the optimal 
value of \ref{LSIP} in a monotonic manner.

An adaptive refinement algorithm is developed to obtain an approximate solution to \ref{LSIP}
which is proved to terminate in finite iterations for arbitrarily given tolerances.
Numerical results show that the algorithm works well in finding feasible solutions for 
\ref{LSIP}.




\begin{thebibliography}{99}

\bibitem{yy1}
S. Christensen, \emph{A method for pricing American options using Semi-infinite linear programming},
Mathematical Finance, 24 (2014), pp.~156--172.


\bibitem{yy2}
S. Daum and R. Werner, \emph{A novel feasible discretization method for linear semi-infinite programming applied to basket option pricing},
Optimization, 60 (2011), pp.~1379--1398.


\bibitem{yy3}
S. \"{O}z\"{o}\u{g}\"{u}r-Aky\"{u}z and G. W. Weber, 
\emph{Infinite kernel learning via infinite and semi-infinite programming},
 Optimisation Methods \& Software, 25 (2010), pp.~937--970.


\bibitem{mma}
M. A. Goberna,and M. A. L\'{o}pez,
\emph{Linear Semi-Infinite Optimization},
Wiley, New York, 1998.

\bibitem{ma} 
M. A. Goberna, and  M. A. L\'{o}pez,
\emph{Linear semi-infinite programming theory: an updated survey}, European Journal of Operational Research,
143 (2002), pp.~390--405.


\bibitem{hkko}
R. Hettich and K. O. Kortanek,
\emph{Semi-infinite programming: theory, methods, and applications},
SIAM Rev., 35 (1993), pp.~380--429.


\bibitem{rg}
R. Reemtsen and S. G\"{o}rner,
\emph{Numerical methods for semi-infinite programming: A survey},
in Semi-Infinite Programming, R. Reemtsen and J.-J. R\"{u}ckmann, eds., Kluwer, Boston,
1998, pp.~195--275.

\bibitem{BB} 
B. Betr\`{o}, \emph{ An accelerated central cutting plane algorithm for linear semi-infinite programming}, Mathematical programming, 101 (2004), pp.~479--495.


\bibitem{sa}
S. {\AA}. Gustafson, \emph{ On the computational solution of a class of generalized moment problems},
SIAM Journal on Numerical Analysis, 7
(1970), pp.~343--357.


\bibitem{tse}
T. Le¨®n, S. Sanmatias and E. Vercher, 
\emph{On the numerical treatment of linearly constrained semi-infinite optimization problems},
European Journal of Operational Research, 121 (2000) pp.~78--91.



\bibitem{ea}
E. J. Anderson and  A. S. Lewis, \emph{An extension of the simplex algorithm for semi-infinite
linear programming}, Mathematical Programming, 44
(1989), pp.~247--269.


\bibitem{em}
E. J. Anderson and M. A. Goberna, \emph{Simplex-like trajectories on quasi-polyhedral sets},
Mathematics of Operations Research, 26
(2001), pp.~147--162.


\bibitem{floudas}
C. A., Floudas, O. Stein \emph{The adaptive convexification algorithm: a feasible point method for semi-infinite programming},
SIAM Journal on Optimization, 18(4), 
(2007) pp.~1187-1208.


\bibitem{wang}
S., Wang, Y.,  Yuan \emph{Feasible method for semi-infinite programs},
SIAM Journal on Optimization, 25(4), 
(2015) pp.~2537-2560.

\bibitem{Amit}
A., Mitsos, Y.,  Yuan \emph{Global optimization of semi-infinite programs via restriction of the right-hand side},
Optimization, 60(10-11)
(2011) pp.~:1291-308.




\bibitem{IT}
G. Alefeld and G. Mayer, \emph{Interval analysis: theory and applications}, J. Comput. Appl. Math.,
121 (2000), pp.~421--464.



\bibitem{postlsip} 
M. A. Goberna,
\emph{Post-optimal analysis of linear semi-infinite programs},
Optimization and Optimal Control, Springer New York, 2010, pp.~23--53.


\bibitem{advance}
M. A. Goberna, \emph{Linear semi-infinite optimization: recent advances}, In Continuous Optimization, Springer US, 2005, pp.~3--22.



\bibitem{KS}
K. Glashoff and S. A. Gustafson, \emph{Linear optimization and approximation},
Springer-Verlag, Berlin, 1983.





\bibitem{TE}
T. Leon and E. Vercher, \emph{A purification algorithm for semi-infinite programming},
European Journal of Operational Research, 57 (1992), pp.~412--420.






\bibitem{IT1}
R. Moore, \emph{Methods and applications of interval analysis},
SIAM, Stud. Appl. Math. 2, Philadelphia, 1979.


\bibitem{convex}
R. T. Rockafellar, \emph{Convex analysis},
Princeton University Press, New Jersey, 1970.

\bibitem{intlab}
S. M. Rump, \emph{INTLAB - INTerval LABoratory, Institute for Reliable Computing, Hamburg
University of Technology}, 1999, http://www.ti3.tu-harburg.de/rump/intlab.


\end{thebibliography}

\section{Appendices}
\noindent\textbf{Proof of Lemma \ref{lm32}}\medskip

Since the Slater condition holds, there exists a point $\bar{x}\in \mathbb{R}^n$ such that
\[a(y)^{\top}\bar{x}+a_0(y) > 0, \ \forall y\in Y.\]
It has been shown in \cite{ma} that the boundary of $F$ is
\[\partial F = \{x\in F\ |\ \min_{y\in Y}\{ a(y)^{\top}\bar{x}+a_0(y)\} = 0 \}.\]
It follows that $F = F^o \cup \partial F$. The compactness of the index set $Y$ implies that the function
$g(x) = \min_{y\in Y}\{ a(y)^{\top}x+a_0(y)\}$ is continuous. Thus, $F$ is closed.
It suffices to prove that
\[\partial F \subseteq cl(F^o).\]
For any $\tilde{x}\in \partial F$, we have $a(y)^{\top}\tilde{x}+a_0(y) = 0$ for all $y\in A(\tilde{x})$ with
$A(\tilde{x}) = \{y\in Y\ |\ a(y)^{\top}\tilde{x}+a_0(y) = 0\}$.
Then
\[a(y)^{\top}(\bar{x}-\tilde{x}) > 0, \ \forall y\in A(\tilde{x}).\]
This indicates that for any $\tau > 0$, we have
\[a(y)^{\top}(\tilde{x} + \tau(\bar{x}-\tilde{x})) + a_0(y) > 0,\ \forall y\in A(\tilde{x}).\]
For a point $y\in Y$ and $y\notin A(\tilde{x})$, there holds that $a(y)^{\top} \tilde{x} + a_0(y) > 0$.
Therefore, $a(y)^{\top}(\tilde{x} + \tau(\bar{x}-\tilde{x})) + a_0(y) > 0$ for $\tau$ small enough.
Since $Y$ is compact, we can chose a uniform $\tau$ such that 
$a(y)^{\top}(\bar{x} + \tau(\bar{x}-\tilde{x})) + a_0(y) > 0, \ \forall y\in Y$ for $\tau$ 
small enough. It follows that we can choose a sequence 
$\tau_k >0$ with $\lim_{k\to\infty} \tau_k = 0$ such that
\[a(y)^{\top}(\tilde{x} + \tau_k(\bar{x}-\tilde{x})) + a_0(y) > 0, \forall y\in Y, k\in \mathbb{N}.\]
Hence $x_k = \tilde{x} + \tau_k(\bar{x}-\tilde{x})\in F^o$ and $\lim_{k\to \infty} x_k = \tilde{x}$ which implies
that $\tilde{x}\in cl(F^o)$.

This completes our proof.\medskip
\newline
\newline

\noindent\textbf{Proof of Lemma \ref{lm33}}\medskip

The Slater condition implies that there exists a point $\bar{x}\in \mathbb{R}^n_{+}$ such that
\[a(y)^{\top}\bar{x} + a_0(y) > 0, \ \forall y\in Y.\]
Let $T = \{\tau_k\ |\ k = 0,1,...,N\}$, from (\ref{haus}) we know that for each $Y_k = [\tau_{k-1},\tau_k]$,
$k = 1,2,...,N$, there holds that
\[\min_{y\in Y_k} a_i(y) - A_{i,k}^l \leq \gamma_i |Y_k|^p \leq \gamma_i |T|^p, i = 0,1,...,n, k = 1,2,...N,\]
with $p\geq1$.
By direction computation, we have
\begin{equation*}
  \sum_{i = 1}^n [\min_{y\in Y_k} a_i(y)]\bar{x}_i + \min_{y\in Y_k} a_0(y) - [\sum_{i = 1}^n A_{i,k}^l \bar{x}_i
  + A_{0,k}^l] \leq [\sum_{i = 1}^n \gamma_i \bar{x}_i + \gamma_0] |Y_k|^p.
\end{equation*}
The Lipschitz continuity of $a_i(y)$, $i = 0,1,...,n$, implies that
\[\min_{y\in Y_k} [\sum_{i = 1}^n a_i(y)\bar{x}_i + a_0(y)]-[\sum_{i = 1}^n [\min_{y\in Y_k} a_i(y)]\bar{x}_i +
\min_{y\in Y_k} a_0(y)] \leq [\sum_{i = 1}^n L_i\bar{x}_i + L_0]|Y_k|.\]
It follows from the last two inequalities that
\begin{equation*}
  \sum_{i = 1}^n A_{i,k}^l \bar{x}_i + A_{0,k}^l \geq \min_{y\in Y_k} [\sum_{i = 1}^n a_i(y)\bar{x}_i + a_0(y)]
  - \{[\sum_{i = 1}^n \gamma_i \bar{x}_i + \gamma_0] |Y_k|^p + [\sum_{i = 1}^n L_i\bar{x}_i + L_0]|Y_k|\},
\end{equation*}
which implies that $\sum_{i = 1}^n A_{i,k}^l \bar{x}_i + A_{0,k}^l \geq 0$, $k = 1,2,...,N$, for $|T|$ small
enough. This implies that $\bar{x}$ is a feasible point for the approximate region $F(T)$.

This completes our proof.\medskip
\newline
\newline

\noindent\textbf{Proof of Theorem \ref{t34}}\medskip

By the construction of the approximate regions, we know that $F(T_k)$ is included in the original feasible
set, i.e., $F(T_k)\subseteq F$ for all $k\in \mathbb{N}$. Hence, we have
$\{x_k^* \} \subseteq F$.

Let $\bar{x}$ be any Slater point, we can conclude from Lemma \ref{lm33} that $\bar{x}$ is contained 
in $F(T_k)$ for $k$ large enough. Thus $c^{\top}x_k^* \leq c^{\top}\bar{x}$ which indicates that 
$x_k^*\in L(\bar{x})$ for sufficient large $k$.
Since the level set $L(\bar{x})$ is compact, there exists at least an accumulation point $x^*$ of
the sequence $\{x_k^*\}$. Assume without loss of generality that the sequence $\{x_k^*\}$ itself 
converges to $x^*$, i.e., $\lim_{k\to \infty} x_k^* = x^*$.
It suffices to prove that $x^*$ is an optimal solution to \ref{LSIP}. It is obvious that $x^*$ is feasible for \ref{LSIP}.

Let $x_{opt}$ be an optimal solution to \ref{LSIP}. If $x_{opt}\in F^o$,
then $x_{opt}\in F(T_k)$ for all $k$ large enough. This indicates that $f(x_k^*) \leq f(x_{opt})$ for $k$
large enough and thus
\[f(x^*) = \lim_{k\to \infty} f(x_k^*) \leq f(x_{opt}),\]
where $f(x) = c^{\top}x$.
If $x_{opt}$ lies on the boundary of the feasible set $F$, there exists a sequence of the Slater points
$\{\bar{x}_j\ |\ \bar{x}_j\in F^o\}$ such that $\lim_{j\to \infty} \bar{x}_j = x_{opt}$. For each
$\bar{x}_j\in F^o$ there exists at least an index $k = k(j)$ such that $\bar{x}_j\in F(T_{k(j)})$
which implies that $f(x_{k(j)}^*) \leq f(\bar{x}_j)$ for $j\in \mathbb{N}$.
Since $\{x_k^*\}$ converges to $x^*$ and $\{x_{k(j)}^*\}$ is a subsequence of $\{x_k^*\}$, the sequence
$\{x_{k(j)}^*\}$ is convergent and $\lim_{j\to \infty} f(x_{k(j)}^*) = f(x^*)$.
By the continuity of $f$ we have
\[f(x^*) = \lim_{j\to \infty} f(x_{k(j)}^*) \leq \lim_{j\to \infty} f(\bar{x}_j) = f(x_{opt}).\]
To sum up, we have $x^*\in F$ and $f(x^*) \leq f(x_{opt})$.

This completes our proof.\medskip
\newline
\newline

\noindent\textbf{Proof of Lemma \ref{l35}}\medskip

Let $\bar{x}\in F $ be a Slater point, then we have 
\[a(y)^{\top}\bar{x} + a_0(y)  > 0, \ \forall y\in Y_k, k = 1,2,...,N.\]
Since $a_i(\cdot), i = 1,2,...,n$ are twice continuously differentiable, they are Lipschitz continuous,
i.e., there exist a constant $L$ such that
\[|a_i(y) - a_i(z)| \leq L|y-z|,\ \forall y,z \in Y.\]
Let $\bar{g}_k(x) = \sum_{i=1}^n \min\{\bar{a}_i(\tau_{k-1}),\bar{a}_i(\tau_k)\} x_i + \min\{\bar{a}_0(\tau_{k-1}),
\bar{a}_0(\tau_k)\}$, then we have
\begin{align*}
& a(y)^{\top}\bar{x} + a_0(y) - \bar{g}_k(\bar{x})  \\
& = \sum_{i=1}^n a_i(y) \bar{x}_i + a_0(y) - 
\sum_{i=1}^n \min\{\bar{a}_i(\tau_{k-1}),\bar{a}_i(\tau_k)\} \bar{x}_i + \min\{\bar{a}_0(\tau_{k-1}),
\bar{a}_0(\tau_k)\} \\
& \leq (L\sum_{i=1}^n (\bar{x}_i +1)) |Y_k|, \ \forall y\in Y_k, k = 1,2,...,N. 
\end{align*}
It follows that $\bar{g}_k(\bar{x}) \geq 0$ if $|T|$ is sufficiently small which implies
$\bar{x} \in \bar{F}(T)$.

This completes our proof. \medskip
\newline
\newline
\noindent\textbf{Convergence of Algorithm 1 for R1-LSIP}\medskip

Since $x_k$ is a solution of R1-LSIP($T_k$) for a consistent subdivision $T_k = \{\tau^k_j\ | \ j = 0,1,...,N_k\}$
in the $k$th iteration, it must satisfy the KKT condition as follows:
\begin{equation}\label{kktapp1}
  c - \sum_{j\in A(x^*_k)} \lambda^k_j A_{T_k}(:,j) = 0
\end{equation}
where $A(x^*_k) = \{j\ |\ A(:,j)^{\top} x^*_k + b_{T_k}(j) = 0\}$ and 
$A_{T_k}(i,j) = A^l_{i,j}$, $b_{T_k}(j) =
A^l_{0,j}$ is the corresponding lower bound for $a_i(y)$ and $a_0(y)$ on $[\tau^k_{j-1},\tau^k_j]$.
By (\ref{haus}) we know that for any $\bar{\tau}^k_j\in [\tau^k_{j-1},\tau^k_j]$ there holds that
\[|a_i(\bar{\tau}^k_j) - A_{T_k}(i,j)| \leq \gamma^k_i |\tau^k_j - \tau^k_{j-1}|^p, 0\leq i\leq n, j\in A(x^*_k)\]
where $\gamma^k_i, 1\leq i\leq n$, $p\geq 1$ are constants. Thus, there exist some constants
$0 \leq \beta^k_i \leq \gamma_i^k, i = 0,1,2,...,n$ such that 
$A_{T_k}(i,j) = a_i(\bar{\tau}^k_j)  + \beta^k_i |\tau^k_j - \tau^k_{j-1}|^p$. Substitute this into
(\ref{kktapp1}) and $A(x^*_k)$ we have
\begin{align*}
& c - \sum_{j\in A(x^*_k)} \lambda^k_j [a(\bar{\tau}^k_j) + (|\tau^k_j - \tau^k_{j-1}|^p)\beta^k] = 0, \\
& A(x^*_k) = \{j\ |\ [a(\bar{\tau}^k_j) + (|\tau^k_j - \tau^k_{j-1}|^p)\beta^k] x^*_k +
a_0(\bar{\tau}^k_j) + (|\tau^k_j - \tau^k_{j-1}|^p)\beta_0^k = 0\},
\end{align*}
where $a(\bar{\tau}^k_j) = (a_1(\bar{\tau}^k_j), a_2(\bar{\tau}^k_j),...,a_n(\bar{\tau}^k_j))^{\top}$.
It follows that $x^*_k$ is a $(\epsilon,\delta)$ KKT point of \ref{LSIP} if
the lengths of the subsets $[\tau^k_{j-1},\tau^k_j]$ for $j\in A(x^*_k)$ converge to zero as $k$
goes to infinity. This is true due to the similar argument in the proof of theorem \ref{conapp2}.

\end{document}